\documentclass{amsart}

\usepackage{amssymb}
\usepackage{amsmath}
\usepackage{amsthm} 

 \newcommand{\HH}{{\mathbb{H}}}
 
 \newcommand{\ZZ}{{\mathbb{Z}}}
 \newcommand{\cc}{{\mathcal{C}}}

 \newcommand{\fg}[1]{\pi_1({#1})}
 \newcommand{\py}[1]{{\pi_Y({#1})}}
 \newcommand{\cover}[1]{\widetilde #1}
 \newcommand{\acover}[1]{\widehat #1}
 \newcommand{\Aut}{{\rm Aut}}
 \newcommand{\diam}{{\rm{diam}}}
 \newcommand{\MCG}{{\mathcal{MCG}}} 
 \newcommand{\bdy}{{\partial}} 
 \newcommand{\dcc}{d_\cc}
 \renewcommand{\ni}{\noindent}

\newcommand{\ML}{{\mathcal{ML}}} 
\newcommand{\PML}{{\mathcal{PML}}}
\newcommand{\UML}{{\mathcal{UML}}}

\theoremstyle{plain}
\newtheorem{theorem}{Theorem}[section]

\newtheorem{lemma}[theorem]{Lemma}
\newtheorem*{maintheorem}{Main Theorem}
\newtheorem*{maincorollary}{Main Corollary}

\theoremstyle{definition}
\newtheorem*{define}{Definition}
\newtheorem*{claim}{Claim}

\newtheorem*{remark}{Remark}

\begin{document}

\title{Big Handlebody Distance Implies Finite Mapping Class Group}

\author{Hossein Namazi}

\date{}

\address{Department of Mathematics, Stony Brook University, Stony Brook, NY 11794}

\email{hossein@math.sunysb.edu}

\begin{abstract}
We show that if $M$ is a closed three manifold with a Heegaard 
splitting with sufficiently big {\em handlebody distance} then the subgroup of the mapping class group of the Heegaard surface, which extend to both handlebodies is finite. As a corollary, this implies that under the same hypothesis, the mapping class group of $M$ is finite.
\end{abstract}

\maketitle

\section{Introduction}
\label{sec:Introduction}

It is well known that any closed orientable 3-manifold is obtained by taking two copies of a handlebody and gluing them along the boundary. Such a decomposition is called a {\em Heegaard splitting}. An outstanding problem in studying 3-manifolds is to obtain information about the manifold from the gluing map. One important feature is that when the gluing map is complex, we expect to get a rigid structure for the manifold. The first result in this direction could be Haken's lemma \cite{Haken} which proves that if the gluing map does not take any meridian to another meridian then the manifold is irreducible. Casson and Gordon \cite{CG} generalized this and proved that if in addition to the above assumption we know that image of a meridian has zero intersection with another meridian then the manifold is Haken. Hempel \cite{Hempel} generalized this and defined a {\em Handlebody distance} for a Heegaard splitting. In this definition the assumption in Haken's lemma can be translated to having positive handlebody distance and the Casson-Gordon's assumption is equivalent to having handlebody distance equal to one. He also conjectured that if this number is bigger than two, then the manifold is hyperbolic. Here, we show that if this number is sufficiently big then there are only finitely many isotopy classes of the automorphisms of the Heegaard surface which can be extended to both handlebodies. Then we use this and by applying a result of Jaco and Rubinstein \cite{JR}, we show that the mapping class group of the manifold is finite.

Let $M$ be a closed orientable 3-manifold with a Heegaard splitting. In other terms, $M=H^+\cup_S H^-$ is the union of two handlebodies $H^+$ and $H^-$ with the same genus, glued along their boundaries. We call $S=\bdy H^+=\bdy H^-$ a {\em Heegaard surface} and it is defined up to isotopy. 

To define the handlebody distance for the splitting, we need to define the {\em curve complex} of $S$, $\cc(S)$. This is a locally infinite simplicial complex, whose vertices are homotopy classes of essential simple closed curves and $k$-simplices are $(k+1)$-tuples $[\alpha_0,\ldots,\alpha_k]$ that have simultaneously disjoint representatives.
It is also equipped with a path metric $\dcc$ on its one-skeleton, which makes the one-skeleton of each simplex regular Euclidean with side-length 1 and is the shortest path distance between any two points of the one-skeleton.

Each handlebody $H^+$ and $H^-$ has a set of meridians as a subset of
$\cc(S)$; these are essential simple closed curves which are compressible 
in the handlebody. We call them $\Delta^+$ and $\Delta^-$ respectively. The {\em handlebody distance} for the splitting is 
$\dcc(\Delta^+,\Delta^-)$, their distance in the curve complex of $S$.

As usual we denote the set of isotopy classes of automorphisms of $S$ by $\MCG(S)$. We also will need to consider the set of isotopy classes of the homeomorphisms of the 3-manifold $M$ to itself, which we call the {\em mapping class group} of $M$ and denote by $\MCG(M)$. We say an automorphism of boundary of a handlebody {\em extends} to the handlebody if it is the restriction of an automorphism of the handlebody to the boundary. It is easy to see that if two automorphisms of the boundary are isotopic and extend to the handlebody then their extensions are isotopic. To see this, take a system of compressing disks which cut the handlebody to a ball. We can see that the extensions to this system are isotopic because the handlebody is irreducible. After doing this isotopy, what remains is a ball and two automorphisms which are isotopic on its boundary, but any two such automorphisms are isotopic in the interior of the ball; so we can extend the isotopy to the ball and therefore to the entire handlebody. 

We are interested in a subgroup of $\MCG(S)$, which consists of those elements whose representatives extend to both $H^+$ and $H^-$. We denote this subgroup by $\Gamma(H^+,H^-)$. The above argument shows that there is a well defined map 
$$\Phi:\Gamma(H^+,H^-)\rightarrow \MCG(M)$$
by extending the automorphism of the surface to both handlebodies.

\begin{maintheorem}
If $M$ is a closed 3-manifold which admits a Heegaard splitting $M=H^+\cup_SH^-$ with handlebody distance bigger than $2K+2\delta$ then $\Gamma(H^+,H^-)$, the subgroup of $\MCG(S)$ whose elements extend to both handlebodies, is finite.  

Here $K$ and $\delta$ depend only on the genus of the splitting and 
will be described in section \ref{sec:Preliminaries}.
\end{maintheorem}

Thompson \cite{Thomp} and Hempel \cite{Hempel} showed that any Heegaard splitting of a toroidal manifold has handlebody distance at most $2$. Therefore, the manifolds satisfying the hypothesis of the above theorem are all atoroidal. On the other hand, in a widely accepted but unpublished result, Jaco and Rubinstein have proved that an atoroidal 3-manifold has at most a finite number of splittings of the same genus up to isotopy \cite{JR}. Assuming this result, we can prove our main application of the above theorem:

\begin{maincorollary}
If a closed 3-manifold admits a Heegaard splitting with sufficiently large handlebody distance then the mapping class group of the manifold is finite. 
\end{maincorollary}

\begin{proof}
An image of a Heegaard splitting by a homeomorphism is a Heegaard splitting and if two homeomorphisms are isotopic, the images will give the same splittings up to isotopy. 
Jaco and Rubinstein's result shows that the subgroup of the mapping class
group of $M$ which preserves the Heegaard surface up to isotopy has finite index. To prove our claim we need to show this subgroup is finite. By possibly taking 
another index two subgroup, we can even restrict to the subgroup of those
mapping classes that preserve the surface up to isotopy and after changing it by an isotopy to make it preserve the surface, it also preserves each handlebody. In other terms, those which preserve the orientation of the surface. We say these elements {\em preserve the splitting}. 

Now look at the map $\Phi:\Gamma(H^+,H^-)\rightarrow \MCG(M)$ defined above. We claim that the image is the subgroup of $\MCG(M)$ of all elements that preserve the splitting. In other terms, we have a surjective map from the subgroup of $\MCG(S)$ that extends to both handlebodies to the subgroup of $\MCG(M)$ that preserves the splitting. Take a representative of an element of $\MCG(M)$ that preserves the splitting. We can assume that it preserves $S$, thus induces an automorphism of the surface and an element of $\MCG(S)$. This of course extends to both handlebodies by definition and therefore is in $\Gamma(H^+, H^-)$. It is clear that image of this by $\Phi$ recovers the element of $\MCG(M)$, which we started with and we have proved the claim.

The theorem above shows that if the handlebody distance is large $\Gamma(H^+, H^-)$ is finite and so is its image by $\Phi$. This proves the finiteness of $\MCG(M)$ by what we said in the beginning.
\end{proof}

\begin{remark}
In the above proof, we showed that $\Phi$ is a surjective map from $\Gamma(H^+,H^-)$ to the subgroup of $\MCG(M)$ whose elements preserve the splitting. In general, this map does not need to be injective. As an example, take any standard Heegaard splitting of $S^3$ with a surface of genus bigger than $1$. It is easy to see that a separating simple closed curve, which is boundary of a disk in one of the handlebodies, will bound a disk in the other one as well. Hence, the Dehn twist about this curve gives an automorphism of the surface that extends to both handleobdies and therefore $\Gamma(H^+, H^-)$ is nontrivial, but $\MCG(S^3)$ is trivial.

One can ask if this kernel is trivial for sufficiently big handlebody distance. We have an argument to show that it is trivial when the handlebody distance is sufficiently big and the splitting has some {\em bounded geometry} condition, which can be described in terms of the combinatorics of $\Delta^+$ and $\Delta^-$ in the curve complex. On the other hand, one may be able to show this kernel is trivial or just finite when the splitting is {\em strongly irreducible} (the handlebody distance is at least $2$).
\end{remark}

We should remark that knowing the Geometrization Conjecture implies Hempel's conjecture, that $M$ is hyperbolic when the handlebody distance for a Heegaard splitting of $M$ is bigger than $2$. When $M$ is hyperbolic, an easy consequence of Mostow rigidity would be to see that the set of homotopy classes of the homeomorphisms of $M$ to itself is finite. This is a quotient of $\MCG(M)$, the way we have defined it. Still, work of Gabai-Meyerhoff-Thurston \cite{GMT} proves that for closed hyperbolic three manifolds $\MCG(M)$ is also finite. Hence, a proof of geometrization conjecture proves the above corollary for handlebody distance bigger than $2$. We should also mention that Lustig and Moriah \cite{LM} have proved a similar result by assuming that the splitting has a property which they call {\em double rectangle condition}.

Note that when $S$ is a torus, it is well known that $M$ is $S^3$, 
$S^2\times S^1$ or a lens space. All these cases are well understood
and therefore we always assume that the genus of the splitting is at least 
two.

I want to thank Yair Minsky for helping me with this project as my dissertation advisor, sharing his enlightening ideas and also laying the foundations on which the research in this paper is built.


\section{Preliminaries}
\label{sec:Preliminaries}

In \cite{MM1}, \cite{MM2} and \cite{MM3}, Masur and
Minsky studied various properties of the curve complex.
The contents of this section are a few of those which we will use
in the course of our proof. 

In the introduction, we defined the curve complex for closed surfaces
of genus at least two. 
We extend the definition for any finite type surface $S_{g,b}$, the surface of genus $g$ with $b$ punctures. In
each case, it is a locally infinite simplicial complex with a path metric.

When $3g+b\geq5$, we take {\em non-peripheral} simple closed curves
and connect two of them with an edge, if they have disjoint representatives. 
Here, non-peripheral curves are essential curves that are not boundary parallel.
For $S_{1,0}$, $S_{1,1}$ and $S_{0,4}$, we again consider the set of 
non-peripheral curves but we connect two vertices $\alpha$ and $\beta$, 
if they intersect minimally, i.e. $i(\alpha,\beta)=1$ for $S_{1,0}$ and 
$S_{1,1}$ and $i(\alpha,\beta)=2$ for $S_{0,4}$. The only other case we need 
is for the closed annulus, for which we give the definition below
and for all other cases we define the curve complex to be empty.

\noindent
{\it Definition of the curve complex for annulus.} 
Consider an annulus $A$ with its two boundaries. 
The vertices of $\cc(A)$ are the homotopy classes of arcs
connecting these two boundaries relative their end points. This of course
will be an uncountable set of vertices; we connect two vertices with
an edge when they have representatives with disjoint interiors. 

The metric $\dcc$ on the one-skeleton of the curve complex is defined as before: the one-skeleton of each simplex is regular Euclidean of side-length 1 and the distance between any other two points in the one-skeleton is the length of the shortest path connecting them. 
For a closed annulus, it is not hard to see that the complex with its metric is {\em quasi isometric} to $\ZZ$ with its word metric \cite{MM2}.
Note that in general, one can extend the metric to get a distance function between any two points of the curve complex, but this will be always quasi-isometric to the one-skeleton with the metric defined above.

From now on all the surfaces that we consider will be $S_{g,b}$, where
$3g+b \geq 4$ and the subsurfaces will be those plus annular subsurfaces. We
also require the subsurfaces to be essential which means that all of their
boundary components has to be non-peripheral.

If $S$ is a surface of finite type, we define its {\em mapping class
group} $\MCG(S)$, to be the set of isotopy classes of orientation
preserving automorphisms of $S$ which fix any component of the boundary.

Thurston's classification of surface automorphisms divides elements of 
$\MCG(S)$ into periodic, reducible and pseudo-Anosov ones; 
where  periodic ones have finite order and reducible ones are those which 
preserve a set of homotopy classes of disjoint 
non-peripheral simple closed curves. The rest are pseudo-Anosovs. If 
$\phi \in \MCG(S)$ is a pseudo-Anosov it has a stable and unstable laminations
$\lambda^+$ and $\lambda^-$; these are two transversal measured laminations 
that fill the surface and for any $\alpha \in \cc(S)$, $\phi^n(\alpha)$ 
converges to $\lambda^+$ in the space of projectivised measured laminations 
$\PML(S)$, when $n\rightarrow \infty$ and converges to $\lambda^-$, when 
$n\rightarrow -\infty$. For an exposition of Thurston's theorem and more on the spaces of laminations see \cite{CB} and \cite{FLP}.

\begin{lemma}[Luo, Masur-Minsky 1999 \cite{MM1}]
\label{luo}
Let $Y$ be any surface, possibly with boundary and $\phi$ any pseudo-Anosov 
acting on $Y$. For any $\alpha\in\cc(S)$ the set $\{\phi^n(\alpha)\}_{n\geq0}$ has infinite diameter in $\cc(Y)$.
\end{lemma}

The next theorem is significant in studying the curve complex. It states 
that $\cc(S)$ defined as above is $\delta$-hyperbolic
in sense of Gromov. This means that in any geodesic triangle, each side is
in a $\delta$ neighborhood of the other two. 

\begin{theorem}[Hyperbolicity, Masur-Minsky 1999 \cite{MM1}]
\label{hyperbolicity}
Let $S$ be a surface of finite type with negative Euler Characteristic and $\cc(S)$ its curve complex. Then $\cc(S)$ with its path metric has infinite diameter and is $\delta$-hyperbolic in sense of Gromov. 
\end{theorem}

\ni
{\it Note:} Here, we always assume that $\delta\geq 2$.

If $Y$ is an essential subsurface of $S$, we can define a projection 
$\pi_Y$ from $\cc_0(S)$ to subsets of $\cc_0(Y)$ of diameter less than four, when $Y$ is non-annular. 
Even when $Y$ is annular we define a projection. For a more precise
definition, see \cite{MM2}.

\begin{define}
If $\alpha \in \cc_0(S)$ does not intersect $Y$ essentially (equivalently if it is isotopic to a curve that does not intersect $Y$ at all)
or $Y$ is a three-holed sphere, we define 
$\py{\alpha}=\emptyset$. If not we have two cases:
\begin{itemize}
\item {\it $Y$ is non-annular.} In this case we consider 
$\alpha \cap Y$. This is a set of disjoint curves and arcs. We replace
arcs by simple closed curves obtained by doing a surgery along $\bdy Y$.
Now we define $\py{\alpha}$ to be the set of all essential simple
closed curves obtained. This is 
nonempty; otherwise, $\alpha \cap Y$ would be inessential. Also, one 
can prove that its diameter in $\cc(Y)$ is at most three.
\item {\it $Y$ is an annulus.} Because of our restriction to 
finite type surfaces with $3g+b\geq4$, $S$ admits
a hyperbolic metric. So, we can identify the universal cover of $S$ with $\HH^2$, which has a compactification as a closed disk and $\fg{S}$ acts on it by isometries. Take the annular cover $\cover Y=\HH^2/\fg{Y}$ of $S$ to which $Y$ lifts homeomorphically. $\fg{Y}$ is a cyclic subgroup of isometries of $\HH^2$ with two fixed points at infinity. The quotient of the closed disk minus these two points is a closed annulus $\acover Y$ that compactifies $\cover Y$ naturally. We consider $\cc(Y)$ to be the same as $\cc(\acover Y)$ and define $\pi_Y$ as a map form $\cc_0(S)$ to set of subsets of $\cc_0(Y)=\cc_0(\acover Y)$ of small diameter. All lifts of the geodesic representative of $\alpha$ to $\cover Y$ naturally obtain arcs in the closed annulus $\acover Y$. We define $\py{\alpha}$ to be the set of those which connect the two boundary components. Again this cannot be empty since $\alpha$ intersects $Y$ essentially and by definition has diameter less than 1.  
\end{itemize}

\ni We also denote the distance between projections of $\alpha$ and $\beta$ in $\cc_0(Y)$
by $d_Y(\alpha,\beta)$ and diameter of projection of a set $D$ to $\cc_0(Y)$ by
$\diam_Y(D)$.
\end{define}

We remark that if $Y$ is a subsurface, any $f\in \MCG(S)$ acts by isomorphism $f:\cc(Y)\rightarrow \cc(f(Y))$, and this fits naturally with the action on $\cc(S)$ via $\pi_{f(Y)}\circ f=f\circ\pi_Y$. In particular, if $Y$ is preserved by $f$ (up to isotopy) then $\pi_Y\circ f= f\circ\pi_Y$, where we have used $f$ for its induced action on $\cc(Y)$ as well. For the annular case, we need the following lemma in parallel to lemma \ref{luo}.

\begin{lemma} \label{dehn twist}
Let $Y$ be annular and contained in a subsurface $Z\subset S$ as a proper and non-peripheral subsurface. Also assume $\phi\in\MCG(S)$ preserves $Z$ and $\phi|_Z$ is a nonzero power of a Dehn twist about the core of $Y$, both up to isotopy. If $\alpha\in\cc_0(S)$ intersects $Y$ essentially then the set $\{\py{\phi^n(\alpha)}\}_{n\geq0}$ has infinite diameter in $\cc(Y)$.
\end{lemma}

\begin{proof}
In \cite{MM2}, it is shown that 
$$d_Y(D_\beta^n(\alpha),\alpha)=2+|n|,$$
for any $\alpha$ that intersects $Y$ essentially, where $\beta$ is the core of $Y$, $D_\beta$ is the positive Dehn twist about $\beta$ and $n\neq0$. This is done after observing the fact that for distinct points $a$ and $b$ in $\cc_0(Y)$, $d_Y(a,b)=1+|a.b|$ where $a.b$ is the algebraic intersection number of their interiors \cite{MM2}. 

Assume $\phi$ is isotopic to $\psi$ a non-zero power of $D_\beta$ in $Z$. The above argument shows that $\{\py{\psi^n(\alpha)}\}_{n\geq0}$ has infinite diameter in $\cc(Y)$. So we just need to show that $\py{\phi^n(\alpha)}$ and $\py{\psi^n(\alpha)}$ have bounded distance from each other. Fix a hyperbolic metric on $S$ and realize all the simple closed curves, including components of $\bdy Z$, by closed geodesics. As before take $\acover Y$ to be the compactification of annular cover of $S$ associated to $Y$ like in the definition of projections. Let $\acover \beta$ be the lift of $\beta$ in $\acover Y$ which is its core and also let $Z'$ be a component of lift of $Z$ that contains $\acover \beta$. Since $Y$ is non-peripheral in $Z$, components of $\bdy Z'$ are properly embedded arcs which do not intersect $\acover \beta$ and their interiors are infinite geodesics which are lifts of $\bdy Z$. Now take a lift of $\phi^n(\alpha)$ to $\acover Y$ that intersects $\acover \beta$ and its compactification $\alpha_n'$, which is a properly embedded arc. Because of the assumption, there exists an isotopy supported on $Z$ between $Z\cap\phi^n(\alpha)$ and $Z\cap\psi^n(\alpha)$ that preserves $\bdy Z$. We can lift this isotopy to $Z'$ such that it preserves $\bdy Z'$. This gives an isotopy in $Z'$ between $Z'\cap\alpha_n'$ and intersection of a lift of $\psi^n(\alpha)$ with $Z'$. Take this lift of $\psi^n(\alpha)$ and we denote its compactification by $\alpha_n''$, which is again a properly embedded arc. $Z'\cap\alpha_n'$ and $Z'\cap\alpha_n''$ are properly embedded arcs which are isotopic in $Z'$ (through properly embedded arcs). Hence, they cannot intersect more than once; otherwise we will have a geodesic bigon between two adjacent intersections, which is impossible. $\alpha_n'$ and $\alpha_n''$ intersect two components of $\acover Y \backslash Z'$; but each of these components is a disk and therefore they cannot intersect more than once in each of them (again because we cannot have geodesic bigons). This shows that overall $\alpha_n'$ and $\alpha_n''$ intersect at most three times and then the identity $d_Y(a,b)=1+|a.b|$ shows that $d_Y(\alpha_n',\alpha_n'')\leq 4$. Then of course, the definition of projections shows that  $d_Y(\phi^n(\alpha),\psi^n(\alpha))\leq 4$ and we are done.
\end{proof}

The next theorem shows the significance of infinite diameter projections. In fact, it shows that if $\alpha$ and $\beta$ are in $\cc(S)$ and have very far projections in $\cc(Y)$, for a proper subsurface $Y\subset S$, then the geodesic connecting them in $\cc(S)$ has distance at most one from $\bdy Y$.

\begin{theorem}[Bounded Geodesic Image, Masur-Minsky 2000 ~\cite{MM2}]
\label{boundedimage}
Let $Y$ be a proper subsurface of $S$ which is not a three punctured sphere and 
let $g$ be a geodesic segment, ray or biinfinite line in $\cc(S)$ such that 
$\py{v} \neq \emptyset$ for every vertex of  $g$.

There is a constant $M$ only depending on the Euler characteristic of $Y$, so that 
$$ \diam_Y(g) \leq M.$$
\end{theorem}

In a geodesic metric space, a {\em $K$-quasiconvex} subset is a subset that
any geodesic connecting two of its points is within distance $K$ of
the subset. 

\begin{theorem}[Quasiconvexity, Masur-Minsky 2003 ~\cite{MM3}] 
\label{quasiconvexity}
If $H$ is a handlebody with boundary $S$ then the set of compressible curves is a 
$K$-quasiconvex subset of $\cc(S)$, where $K$ depends only on the genus of $S$.
\end{theorem}


\section{Proof of the Main Theorem}
\label{sec:ProofOfTheMainTheorem}

\begin{proof}
From now, we assume the handlebody distance is bigger than $2K+2\delta$ where $K$ and $\delta$ are the constants obtained respectively in the quasiconvexity theorem \ref{quasiconvexity} and in the hyperbolicity theorem \ref{hyperbolicity}.
Take $\phi\in \MCG(S)$ that extends to both handlebodies. Since it extends, it has to preserve the set of compressible curves for each handlebody. Thus, if we consider it as an isometry of the curve complex of $S$, it preserves the sub-complexes $\Delta^+$ and $\Delta^-$. This is all we will use to show that the set of such maps is finite.
Using Thurston's classification of elements of $\MCG (S)$, $\phi$ is either 
periodic, reducible or pseudo-Anosov. 

\begin{claim}
If $\dcc(\Delta^+,\Delta^-)$ is sufficiently large and $\phi$ preserves both $\Delta^+$ and $\Delta^-$, then it has to be periodic.
\end{claim}

\begin{proof} We will exclude the other cases.

\vspace{.2cm}
\ni
{\it $\phi$ is not pseudo-Anosov.} Assume it is a pseudo-Anosov. Consider $\alpha \in \Delta^+$ and $\beta \in \Delta^-$ such that $d_\cc (\alpha , \beta )$ is minimum: $d_\cc (\Delta^+,\Delta^-)$. Also consider a geodesic segment $l$ between $\alpha$ and $\beta$ of length $d_\cc(\alpha,\beta)$. Now if we look at image of $\alpha$, $\beta$ and $l$ under $\phi^n$, they converge to $\lambda^+$, the stable lamination of $\phi$, in $\PML(S)$ as $n$ goes to infinity. 

Klarreich's theorem \cite{Klarreich} states that the boundary of the curve complex in sense of Gromov is naturally identified with the set of filling laminations in the set of unmeasured measured laminations of $S$: $\UML(S)$. ($\UML(S)$ is the quotient of $\ML(S)$, where two measured laminations are identified, whenever they have the same support.) We denote this set by $\bdy \cc(S)$. Indeed, she proved that if a sequence of simple closed curves converges in $\PML$ to a projectivised measured lamination supported on $\mu\in\bdy \cc(S)$, it also converges to $\mu$ in $\cc(S)\cup \bdy \cc(S)$ in sense of Gromov. 

Since $\lambda^+$ is filling it projects to an element $[\lambda^+]\in\bdy\cc(S)$ and by what we said this implies that $\phi^n(l)$ converges to $[\lambda^+]$ in $\cc(S)\cup\bdy\cc(S)$ and in particular gets further and further from $l$. Now if we look at the rectangle with vertices $\alpha$, $\beta$, $\phi^n(\alpha)$ and $\phi^n(\beta)$, it has two sides $[\alpha, \beta]$ and $[\phi^n(\alpha), \phi^n(\beta)]$ with fixed length $d_\cc(\Delta^+,\Delta^-)$ and the lengths of the other two sides go to infinity as $n$ goes to infinity. Using hyperbolicity theorem \ref{hyperbolicity}, it is easy to see that in this situation the longer sides get $2\delta$ close. 

On the other hand, $\alpha$ and $\phi^n(\alpha)$ are both in $\Delta^+$ ($\Delta^+$ is invariant under $\phi$), and because $\Delta^+$ is $K$-quasiconvex (Theorem \ref{quasiconvexity}) the geodesic connecting them is within at most $K$ from $\Delta^+$. For the same reason, the opposite side, $[\beta, \phi^n(\beta)]$, is also within $K$ from $\Delta^-$. So $[\alpha,\phi^n(\alpha)]$ is in the $K$-neighborhood of $\Delta^+$ and $[\beta, \phi^n(\beta)]$ is in the $K$-neighborhood of $\Delta^-$ and we showed that they are $2\delta$ close, hence $\Delta^+$ and $\Delta^-$ cannot be more than $2K+2\delta$ apart. This contradicts the assumption about the handlebody distance.

\vspace{.2cm}
\noindent
{\it $\phi$ is not reducible unless it is periodic.} 
If $\phi$ is reducible then we have a system of essential simple closed curves which is invariant by $\phi$. By possibly taking a nonzero power of $\phi$ we may be able to increase the number of such curves. Let $\Lambda = \{ \gamma_1, \gamma_2, \dots, \gamma_k \} $ be a maximal set of non-parallel, pairwise disjoint simple closed curves that is invariant by a power of $\phi$. Such a maximal set exists because we cannot have more than $3g-3$ non-parallel essential simple closed curves on a closed surface of genus $g$.
The idea is to show that $\Lambda$ is within at most $K+1$ both from $\Delta^+$ and $\Delta^-$ and then since $\diam_\cc (\Lambda) \leq 1$, we have
\begin{align*}
d_\cc(\Delta^+,\Delta^-) &\leq d_\cc(\Delta^+,\Lambda) + \diam_\cc (\Lambda) + d_\cc(\Lambda, \Delta^-) \\
& \leq (K+1)+1+(K+1) = 2K+3 \\
& \leq 2K+2\delta,
\end{align*} 
which will be contradicting our assumption.

At first, without loss of generality we can replace $\phi$ with a power $\phi^n$ that 
fixes each element of $\Lambda$ and each component of $S-\Lambda$. By a {\em component domain} for $\Lambda$, we mean
either a component of $S-\Lambda$ or an annular neighborhood of one of elements of $\Lambda$. 
One can see that $\phi$, restricted to each non-annular component domain, gives an element of mapping class group of that component and an isometry of its curve complex. In case of annular domains, we need to look at lift of $\phi$ to the corresponding annular cover and its compactification described in definition of projections to get an isometry of the domain's curve complex.
To get a contradiction in case $\phi$ is not periodic, we claim the following: 

\begin{lemma} \label{unbounded projection}
There exists a component domain $Y$ of $\Lambda$ for which $\py{\Delta^+}$ has unbounded diameter in $\cc(Y)$, except when $\phi$ is periodic. 
\end{lemma}

Let's assume this lemma is true and take $\alpha$ and $\beta$ in $\Delta^+$ such that 
$d_Y(\alpha, \beta) \geq M+1$, where $M$ is the bound from theorem \ref{boundedimage}. 
If $g$ is a geodesic connecting these two in $\cc(S)$, by Quasiconvexity Theorem (theorem \ref{quasiconvexity}), it is within $K$ of $\Delta^+$:
$$\dcc(\Delta^+,v)\leq K, \qquad {\rm any}\ v\in g.$$
On the other hand, $\diam_Y(g) \geq M+1$ since $d_Y(\alpha, \beta) \geq M+1$. Using The Bounded Geodesic Image Theorem (theorem ~\ref{boundedimage}), it means that there has to be a vertex $v$ of $g$ for which $\py{v} = \emptyset$. In other terms, $d_\cc (v, \bdy Y) \leq 1$ and therefore $d_\cc(v, \Lambda)\leq 1$. By adding these facts, we have: 
\begin{align*}
d_\cc(\Delta^+, \Lambda) &\leq d_\cc(\Delta^+, v) + d_\cc(v, \Lambda)\\
                   &\leq K + 1
\end{align*}
as we wanted. The same works for $\Delta^-$ and we will get the contradiction.

Now to finish the proof of the claim, we just need to prove lemma \ref{unbounded projection}.

\begin{proof}[Proof of lemma \ref{unbounded projection}.]
Consider $\Lambda$; if it is not a pants decomposition, then there exists a component
$Y$ of $S-\Lambda$ which is not a three-holed sphere. $\phi|_Y$ is an element of the mapping class group of $Y$. We can again use the classification of surface automorphisms and consider different cases where $\phi|_Y$ is periodic, reducible
or pseudo-Anosov. Because of maximality of $\Lambda$, $\phi|_Y$ cannot be periodic
or reducible: in either case, there would be a non-peripheral simple closed curve which
is preserved by a power of $\phi$ and we could add it to $\Lambda$ and the new set would be still preserved by a power of $\phi$. Therefore it is a pseudo-Anosov. Take an element $\alpha \in \Delta^+$ that intersects $Y$ essentially. Then $\py{\alpha}$ is essential in $\cc(Y)$ and 
$(\phi|_Y)^n(\py{\alpha})=\py{\phi^n(\alpha)}$ gives an unbounded subset of $\cc(Y)$, by lemma \ref{luo} and the remark after the definition of $\pi_Y$. Since $\phi$ preserves $\Delta^+$, we have proved the lemma in this case.

The remaining case is when $\Lambda$ is a pants decomposition and $\phi$ fixes each pants and each of its boundaries. Since the action of $\phi$ on each pants is trivial (isotopic to identity), we can characterize the action just by looking at the actions on annular neighborhoods of the curves in $\Lambda$. For doing this, take any $\gamma\in\Lambda$ and cut the surface along all the elements of $\Lambda$ except $\gamma$. The component that contains $\gamma$ is a new subsurface of $S$, which is either a 4-holed sphere or a 1-holed torus. The restriction of $\phi$ to this subsurface is either a Dehn twist along $\gamma$ or is identity up to isotopy. If the latter happens for all curves in $\Lambda$, it is easy to show that $\phi$ is isotopic to identity and we are done. If the former happens for $\gamma$,  take any element $\alpha\in\Delta^+$ that intersects $\gamma$ essentially. Then lemma \ref{dehn twist} shows that the projections $\py{\phi^n(\alpha)}$ have unbounded diameter in $\cc(Y)$, where $Y$ is the annular neighborhood of $\gamma$ and this finishes proof of the lemma. 
\end{proof}

We have finished proof of the claim and $\phi$ has to be periodic.
\end{proof}

Once we know that all the elements are finite, we can use a well known theorem about mapping class group of surfaces.
For more on this and a proof, see \cite{Ivanov}.

\begin{theorem}[Serre]\cite{Serre}
For any $m \geq 3$, the finite index subgroup 
$$ \ker (\MCG (S) \rightarrow \Aut (H_1(S, \ZZ_m)))$$
of $\MCG(S)$ is torsion free.
\end{theorem}

As a simple corollary, any subgroup of $\MCG(S)$ has a finite index subgroup which is torsion free, say its intersection with one of the subgroups in the theorem. This shows that any subgroup of $\MCG (S)$ that all of its elements are torsion has to be finite and we have proved the main theorem.
\end{proof}

\end{document}